\documentclass[12pt]{article}

\title{A system of grabbing particles \\
related to Galton-Watson trees} 
\author{Jean Bertoin\thanks{Laboratoire de Probabilit\'es, 
UPMC, 175 rue du Chevaleret, 75013 Paris; and DMA, ENS Paris, France. Email:
jean.bertoin@upmc.fr} \and 
Vladas Sidoravicius
\thanks{ IMPA,
Estr. Dona Castorina 110, 
Rio de Janeiro, Brasil; and CWI, Kruislaan 413, 1098 SJ
P.O. Box 94079, 1090 GB Amsterdam
The Netherlands. 
Email: vladas@impa.br} 
 \and 
 Maria Eulalia Vares
 \thanks{CBPF. Rua Xavier Sigaud, 150. 22290-180 Rio de Janeiro, RJ. Brasil.
 Email:  eulalia@cbpf.br} }

\date{}

\usepackage{amsfonts,amsmath,amssymb,bbm}
\usepackage{setspace} 
\def\proof{\noindent{\bf Proof:}\hskip10pt}        
\def\QED{\hfill $\Box$}

\font\tenmath=msbm10 scaled 1200
\font\sevenmath=msbm7 scaled 1200
\font\Phiivemath=msbm5 scaled 1200
\newfam\mathfam \textfont\mathfam=\tenmath
\scriptfont\mathfam=\sevenmath \scriptscriptfont\mathfam=\Phiivemath

\def \\ { \cr }

\def \1{1 \mkern -6mu 1} 
\def\N{\mathbb{N}}

\def\P{\mathbb{P}}
\def\PGW{{\mathbb G \mathbb W}^{(\mu)}}
\def\Z{\mathbb{Z}}

\def\T{\mathbb{T}}

\def \+{^{\rm (+)}}

\newtheorem{theorem}{Theorem}

\newtheorem{lemma}{Lemma}

\begin{document}

\maketitle

\begin{abstract}
 We consider  a system of particles with arms that are activated randomly to grab other particles as a toy model for polymerization.
We assume that the following two rules are fulfilled:
once a particle has been grabbed then it cannot be grabbed again, and an arm cannot grab a particle that belongs to its own cluster. We are interested in the shape   of a  typical polymer in the situation when  the  initial number of monomers is large
and the numbers of arms of monomers are given by i.i.d. random variables. Our main result is a  limit theorem for the empirical distribution of polymers, where limit is expressed in terms of a Galton-Watson tree.
  \end{abstract}

\begin{section}{Introduction}

We consider a system of interacting particles which might be used as a toy model for the formation of certain polymers.
Roughly speaking, each particle possesses a certain number of {\it arms} that are activated successively at random to grab other particles.
When an arm becomes active, it creates an oriented edge from the particle with the activated  arm to the particle which has been grabbed, 
so the system describes the evolution of clusters of vertices with some inactive arms
which are connected by oriented edges. As time passes, the number of inactive arms diminishes whereas the connectivity  in the system increases by the creation of new  edges. We impose the rules that once a particle has been grabbed, then it cannot be grabbed again, and  that each time an arm is activated, it grabs a particle which does not belong to its own cluster.
The latter requirement 
impedes the formation of cycles; thus the structure that connects each cluster of particles is  given by a tree and the terminal state obtained when all arms have been activated can be described as a random forest.

This grabbing particle system bears some resemblance with  the so-called random configuration model,
which has been introduced independently by Bollob\'as \cite{Bollo} and Wormald \cite{Worm}; see also Chapter 3 in Durrett \cite{Dur} for further references. The latter is a natural stochastic algorithm which aims at producing a random graph on a set of vertices with pre-described degrees (in general the resulting graph is not simple and may contain loops and multiple edges). Typically, a certain number of stubs is appended to each vertex, and one joins pairwise the stubs uniformly at random to create edges between vertices. Similarities and differences between the two models will be discussed at the end of this work.

We are mainly interested in the case when the number of particles is large, and more precisely, in the situation when 
the empirical distribution of the number of arms converges to some probability measure $\mu$ on $\Z_+$ as the number of particles goes to infinity.
Thinking of these particles as monomers, our goal is to get information about the statistics of  a typical polymer structure when the polymerization process is completed.  
In this direction, we shall consider the slightly different setting where
the sequence of the number of arms is given by  i.i.d. random variables with law $\mu$,
which is simpler to analyze.  Loosely speaking, we will show that when the number of particles goes to infinity, a typical tree (i.e. a tree chosen uniformly at random in the terminal forest) tends to be distributed as a Galton-Watson tree with reproduction law given
by the distribution of the number of arms of a typical particle. This result is perhaps quite intuitive; however providing a rigorous argument is not so straightforward, and this motivates the present work.

The remainder of this paper is organized as follows. The grabbing particle system
is presented formally in the next section, in the situation when the number of particles is finite and the number of arms of each particle is deterministic. In Section 3, we consider the case when the number of arms is random and given by a sequence of i.i.d. variables, and we make the connexion with Galton-Watson processes. More precisely, 
we shall show that conditionally on the total number of arms, the distribution of the terminal state is that of a Galton-Watson forest with a given number of trees, conditioned by its total size. Finally Section 4 is devoted to the asymptotic study of the empirical measure of trees in the terminal state when the number of particles tends to infinity.

\end{section}

\begin{section} {The grabbing particle system}
Consider $n\geq 2$ particles which are labeled by $1,\ldots, n$. Each particle, say $i$, consists in a vertex to which a certain number $x_i\in\Z_+$ of  arms are attached. Arms can be either active or inactive.  An active arm attached to $i$ is an oriented edge $i\to j$ linking $i$ to another vertex $j\neq i$, whereas an inactive arm is an incomplete oriented edge $i\to $  where only the left-end is specified. Initially, every arm is inactive, and once an arm has been activated, it remains active forever.
We assume that the total number of arms is less than the number of particles, i.e.
\begin{equation}\label{E1}
x_1+\cdots + x_n = n-k\qquad \hbox{ for some integer } 1\leq k \leq n.
\end{equation}
We then enumerate these $n-k$ arms uniformly at random, which specifies the order at  which they are activated. 

Specifically, suppose that the first arm to be activated belongs to the particle $i$. Then at time $1$, this arm grabs a particle $j\neq i$ chosen uniformly
at random amongst the other $n-1$ particles. The system at time $1$ thus consists in
$n-2$ isolated particles with inactive arms, and a cluster of two particles, $i$ and $j$, connected by an oriented edge $i\to j$ which is labeled by $1$.
The particle $i$ has still $x_i-1$ inactive arms, and the particle $j$ has kept its $x_j$ inactive arms; inactive arms are labeled by $2, \ldots, n-k$.
We iterate in an obvious independent manner, with the rules that 
every particle can be grabbed by at most one arm, and that each time an arm is activated, it grabs a particle from a different cluster.

 It is convenient to depict the evolution of the system in the upper half-plane by introducing some further orderings  which of course do not alter the dynamics. Specifically,
 imagine that initially, particles lie on the horizontal axis, such that the initial order from the left to the right  is uniformly random. We may then think of the $x_i$ arms attached to the particle $i$ as  vectors 
 making an angle
 $\ell\pi/(x_i+1)$ for $\ell =1,\ldots, x_i$ with the (directed) horizontal axis.
If, as above, the first arm to be activated belongs to the particle $i$
and grabs particle $j$, then at time $1$, the particle $j$ (together with its arms) is shifted
at the right-end of the first arm. See Figure 1 for this first step of the evolution of the system (with $n=7$, $k=2$, $i=6$ and $j=5$).

\begin{picture}(300,120)
\put(15,20){\circle{15}}
\put(40,20){\circle{15}}
\put(65,20){\circle{15}}
\put(90,20){\circle{15}}
\put(115,20){\circle{15}}
\put(140,20){\circle{15}}
\put(165,20){\circle{15}}

\put(15,20){\makebox(0,0){$3$}}
\put(40,20){\makebox(0,0){$5$}}
\put(65,20){\makebox(0,0){$2$}}
\put(90,20){\makebox(0,0){$6$}}
\put(115,20){\makebox(0,0){$1$}}
\put(140,20){\makebox(0,0){$7$}}
\put(165,20){\makebox(0,0){$4$}}

\put(40,27){\vector(0,1){30}}
\put(140,27){\vector(0,1){30}}
\put(165,27){\vector(0,1){30}}
\put(90,27){\vector(1,1){30}}
\put(90,27){\vector(-1,1){30}}

\put(35,40){\makebox(0,0){\small 3}}
\put(70,40){\makebox(0,0){\small 1}}
\put(97,40){\makebox(0,0){\small 4}}
\put(135,40){\makebox(0,0){\small 2}}
\put(160,40){\makebox(0,0){\small 5}}

\put(300,20){\circle{15}}
\put(325,20){\circle{15}}
\put(350,20){\circle{15}}
\put(375,20){\circle{15}}
\put(400,20){\circle{15}}
\put(425,20){\circle{15}}

\put(300,20){\makebox(0,0){$3$}}
\put(325,20){\makebox(0,0){$2$}}
\put(350,20){\makebox(0,0){$6$}}
\put(375,20){\makebox(0,0){$1$}}
\put(400,20){\makebox(0,0){$7$}}
\put(425,20){\makebox(0,0){$4$}}

\put(400,27){\vector(0,1){30}}
\put(425,27){\vector(0,1){30}}
\put(350,27){\vector(1,1){30}}
\put(350,27){\vector(-1,1){30}}
\put(317,63){\circle{15}}
\put(317,63){\makebox(0,0){$5$}}
\put(317,70){\vector(0,1){30}}

\put(312,83){\makebox(0,0){\small 3}}
\put(330,40){\makebox(0,0){\small 1}}
\put(354,40){\makebox(0,0){\small 4}}
\put(391,40){\makebox(0,0){\small 2}}
\put(420,40){\makebox(0,0){\small 5}}

\end{picture}

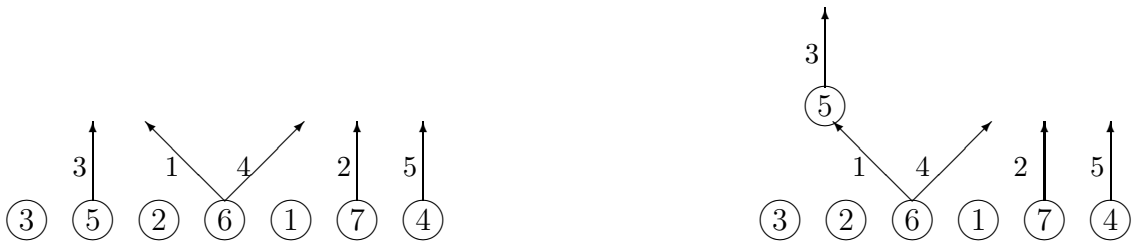
\centerline{\it Figure 1:  planar representation of a grabbing particle system} 
\centerline{\it at times $t=0$ (left) and $t=1$ (right)}
\vskip 4mm

We then continue by an obvious iteration. Only particles lying on the horizontal axis can be grabbed, and  each time such a particle is grabbed by the activation of an arm, it is shifted together with its connected component
to become the right-end of that arm.
As time passes, we thus obtain a growing family of planar rooted trees, where  vertices may bear some further inactive arms, and the roots are the vertices still  lying on the horizontal axis, i.e. those which have not been grabbed so far.
Edges and inactive arms are labeled by 
$1, \ldots, n-k$, the label of an edge being always smaller than that of an inactive arm,
and the vertices are labeled by $1, \ldots, n$.
The evolution reaches its final state at time $n-k$ when all the arms have been activated;
 obviously we then have  a forest of $k$ planar rooted trees with $n-k$ labeled edges on a set of $n$ vertices, and the trees are ordered from the left to the right. See Figure 2 below. We stress that the terminal state retains all the essential\footnote{However the exact ordering of particles lying on the horizontal axis is lost in general, but this is not relevant as far as empirical distributions of the clusters and their structures are concerned.}
  information about evolution of the system, as it suffices to prune the last $n-k-\ell$ edges to recover
 the clusters of particles in the system and their structures at time $\ell$.

\begin{picture}(300,180)(-10,0)
\put(160,20){\circle{15}}
\put(160,64){\circle{15}}
\put(300,20){\circle{15}}
\put(129,109){\circle{15}}
\put(191,109){\circle{15}}
\put(300,64){\circle{15}}
\put(160,71){\vector(1,1){29}}
\put(160,71){\vector(-1,1){29}}
\put(300,27){\vector(0,1){29}}
\put(160,27){\vector(0,1){30}}

\put(129,109){\circle{15}}
\put(129,117){\vector(0,1){30}}
\put(129,154){\circle{15}}

\put(160,63){\makebox(0,0){$6$}}
\put(300,20){\makebox(0,0){$4$}}
\put(129,154){\makebox(0,0){$1$}}
\put(160,20){\makebox(0,0){$7$}}
\put(191,110){\makebox(0,0){$2$}}
\put(300,64){\makebox(0,0){$3$}}
\put(129,110){\makebox(0,0){$5$}}

\put(140,83){\makebox(0,0){\small 1}}
\put(165,83){\makebox(0,0){\small 4}}
\put(295,38){\makebox(0,0){\small 5}}
\put(125,130){\makebox(0,0){\small 3}}
\put(155,38){\makebox(0,0){\small 2}}

\end{picture}

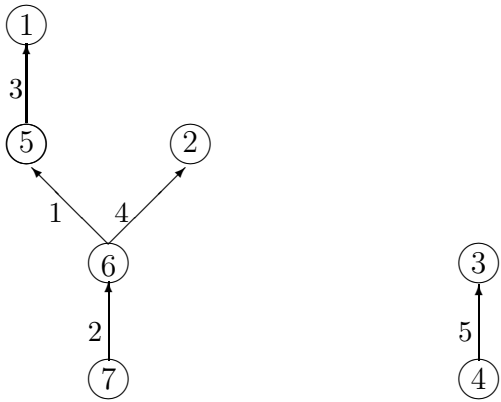
 \centerline{\it Figure 2:   terminal state as a planar forest }
 
 \vskip 4mm
 
We shall now describe the distribution of the terminal state for the grabbing particle system.
In this direction, it is convenient to introduce first the space
$\Phi_{n,k}$ of planar forests with $k$ rooted trees on a set of $n$ labeled vertices, with labeled edges. 
The $k$ roots lie on the horizontal axis and are thus naturally ordered from the left to the right.
The choice of a root in a tree induces an orientation of the edges in that tree, by deciding
that each edge points away from the root.
In this setting, recall the notation $i\to j$ for $i,j\in\{1,\ldots, n\}$ with $i\neq j$ to denote an oriented edge from the vertex labeled by $i$ to the vertex labeled by $j$. For any forest
$\varphi \in \Phi_{n,k}$ and $i\in\{1,\ldots, n\}$, we write
$$d_i(\varphi):={\tt Card}\{j\in\{1,\ldots, n\}\backslash\{i\}: i\to j \hbox{ is an edge of }\varphi\}$$
for  the outer-degree of the
vertex labeled by $i$, i.e. the number of oriented edges which left-end is the $i$-th vertex.

Recall that $x_1, \ldots, x_n$ is a sequence of nonnegative integers such that  \eqref{E1} holds.
We then denote by $\Phi(x_1,\ldots, x_n)\subset \Phi_{n,k}$  the subset of planar forests $\varphi$ with $k$ rooted trees, labeled edges and labeled vertices, such that  the outer-degree  of  $i$ is $d_i(\varphi)=x_i$ for every $i=1,\ldots,n$.  

\begin{lemma}\label{L1} For every  integer $1\leq k\leq n$ and every  sequence of nonnegative integers $x_1, \ldots, x_n$  which fulfills \eqref{E1},
the distribution of the terminal state
of the grabbing particle system  is given by the uniform law on $\Phi(x_1,\ldots, x_n)$.
\end{lemma}

\proof  For $\ell=0,\ldots, n-k$, let $\Phi_\ell(x_1,\ldots, x_n)$ denote the set of planar configurations which can be obtained from the grabbing particle system at time $\ell$, when initially the particle $i$ possesses $x_i$ inactive arms for $i=1,\ldots, n$. In other words,  $\Phi_\ell(x_1,\ldots, x_n)$ is the set of planar forests on a set of $n$ labeled vertices, with $n-\ell$ rooted trees lying on the horizontal axis, to which a total of $n-k-\ell$ yet incomplete edges (i.e. with only a left-end) have been attached to certain vertices, in such a way that for every 
$i=1,\ldots, n$, the number of complete or incomplete edges with left-end $i$ is $x_i$. Further, 
complete or incomplete edges are enumerated in such a way that complete edges are listed before the incomplete ones.
In particular for $\ell=n-k$, in the notation above, there is the identity $\Phi_{n-k}(x_1,\ldots, x_n)=\Phi(x_1,\ldots, x_n)$.

Plainly, under the assumptions made in the claim, the initial state has the uniform distribution in $\Phi_0(x_1,\ldots, x_n)$. We claim that the activation of the first arm yields a random state in $\Phi_{1}(x_1, \ldots, x_n)$ which has again the uniform distribution. Indeed, pick an arbitrary state $\varphi\in \Phi_{1}(x_1, \ldots, x_n)$.
There are exactly $n$ states in $\Phi_{0}(x_1, \ldots, x_n)$ which may possibly yield $\varphi$
by activating the first arm. More precisely, these are the states which are obtained from
$\varphi$ by pruning its unique complete edge at its right-end, and inserting this right-end
(together with its inactive arms)  on the horizontal axis relatively to the $n-1$ roots of $\varphi$.
Further, starting from any of these $n$ states, the probability of getting $\varphi$ after activating their first arm is exactly the probability
that this arm picks the right root amongst the $n-1$ possible ones, that is $1/(n-1)$. 
This establishes our claim. 

We then can iterate this argument, showing that the terminal state is uniformly
 distributed on $\Phi_{n-k}(x_1, \ldots, x_n)= \Phi(x_1, \ldots, x_n)$, as stated. \QED

\end{section}

\begin{section} {Connexion with Galton-Watson forests}

We now turn our attention to the situation when the number of arms of particles is 
given by an i.i.d. sequence of random variables.
Specifically, we consider  a probability measure $\mu$ on $\Z_+$ and  assume that the number of arms of the $i$-th particle, which it is henceforth convenient to denote by $\xi_i$,   is  random with law $\mu$, and that $\xi_1, \ldots, \xi_n$ are independent. 
We shall see that the terminal state of the system can then be described in terms 
of a Galton-Watson forest; in this direction, we first develop some material in this field.

Consider a Galton-Watson process with reproduction law $\mu$, so, roughly speaking,
we have a population model with discrete non-overlapping generations, in which each 
individual gives birth at the next generation to a random number of children according to the law $\mu$ and independently of the other individuals. We shall work on the event  that the total descent of each individual is finite a.s. 
By ranking the ancestors and the progeny of each individual, we can represent
the process by a planar forest with rooted trees, where roots correspond to the ancestors at the initial generation, and each tree describes a genealogy, see e.g. Section 6.2 in Pitman \cite{PiSF}. 

We denote by $F_{n,k}$ the finite set of all planar forests with  $k$ rooted trees
on a set of $n$ vertices (vertices and edges are not yet labeled).
Recall that vertices and edges in a forest  in $F_{n,k}$ can be listed by using, for instance, the depth-first search algorithm (cf. Figure 6.1 in \cite{PiSF} on its page 126).
In this setting, assigning labels to the $n$ vertices (respectively, to the $n-k$ edges)
is equivalent to choosing a permutation on $\{1, \ldots, n\}$ (respectively, a permutation on $\{1, \ldots, n-k\}$). Thus, if  for every $\ell\in\N$,  $\Sigma_{\ell}$ stands for the group of permutations on $\{1,\ldots, \ell\}$, then there is the  canonical identification as the product space
$$\Phi_{n,k}=F_{n,k}\times \Sigma_n\times \Sigma_{n-k}\,.$$ 

We then introduce a product probability measure  on $F_{n,k}\times \Sigma_n\times \Sigma_{n-k}$, denoted by $\PGW_{k,n}$, that can thus be viewed as the law of 
 a triplet of independent variables. The first random variable
with values in $ F_{n,k}$ is a Galton-Watson forest  with reproduction law $\mu$, having $k$ trees and conditioned to have a total size $n$
(assuming implicitly that this event has positive probability), and the second and third
are simply uniformly distributed on $\Sigma_n$ and $\Sigma_{n-k}$, respectively.
We may thus think of  $\PGW_{k,n}$ as the law on $\Phi_{n,k}$ 
of a random state which is obtained by picking a planar
forest with $k$ rooted trees on a set of $n$ vertices according to a Galton-Watson process with reproduction law $\mu$, started with $k$ ancestors and conditioned
to have total size $n$, and independently enumerating the edges and the vertices uniformly at random. 

\begin{theorem}\label{T1} Let $\xi_1, \ldots$ be a sequence of i.i.d. variables with law $\mu$. Fix integers $1\leq k\leq n$ and work conditionally on 
$\xi_1+\cdots+ \xi_n=n-k$, provided that the probability of this event is nonzero. Then the distribution of the terminal state 
of the grabbing particle system started with $n$ particles and such that the initial
number of inactive arms of the $i$-th particle is $\xi_i$ for every $i\in\{1,\ldots, n\}$,
is the law $\PGW_{n,k}$ on $\Phi_{n,k}$.

\end{theorem}

\proof  We first work with a Galton-Watson process with reproduction law $\mu$ started with $k$ ancestors; the notation  $\P_k$ will refer to the  distribution of the process.
Let $(X_i)_{i\geq 0}$ be a sequence of i.i.d. variables with law $\mu$,
introduce the downwards skip-free random walk $S_\ell=X_1+\cdots + X_{\ell}-\ell$
and the first-passage time 
$T:=\min\{\ell\geq 1: S_{\ell}=-k\}$.
By listing individuals in the Galton-Watson process according to the depth-first search algorithm, it is well-known that $\P_k$ can be identified with the distribution of the
stopped sequence $(X_i)_{1\leq i \leq T}$. More precisely, in this setting, $X_i$ represents the number of children of the $i$-th individual listed by depth-first search.

Next pick an arbitrary planar forest $f\in F_{n,k}$ with $k$ rooted trees and $n$ vertices
(vertices and edges are not labeled).
The depth-first search algorithm enables us to encode $f$ by some sequence of integers
$(y_i)_{1\leq i\leq n}$
(more precisely, $y_i$ is the outer-degree of the $i$-th vertex found by the depth-first search algorithm), and the probability that the genealogical forest induced by the 
$k$ ancestors of the Galton-Watson process is given by $f$ thus equals
$\mu(y_1) \cdots  \mu(y_n)$. As a consequence, for every permutations $\sigma\in\Sigma_n$ and $\varsigma\in\Sigma_{n-k}$, we have
$$\PGW_{n,k}\left((f,\sigma,\varsigma)\right)= \frac{\mu(y_1) \cdots  \mu(y_n)}{n! (n-k)! z}$$
with $z=\P_k(T=n)$.

Then pick a planar forest $\varphi$ with $k$ rooted trees on a set of $n$ labeled vertices, with labeled edges, which belongs to the subset $\Phi(x_1,\ldots, x_n)$. We can identify
$$\varphi = (f, \sigma, \varsigma)$$
for a unique $f\in F_{n,k}$, $\sigma\in\Sigma_n$ and $\varsigma\in\Sigma_{n-k}$.
More precisely, $f$ is encoded via depth-first search by the sequence 
$(y_i)_{1\leq i\leq n}$ where
$x_i=y_{\sigma(i)}$, and therefore
$$\PGW_{n,k}(\varphi)= \frac{\mu(x_1) \cdots  \mu(x_n)}{n! (n-k)! z}.$$

This shows that, provided that $\mu(x_i)>0$ for every $i=1,\ldots, n$,  the conditional probability $\PGW_{n,k}(\cdot \mid \Phi(x_1,\ldots, x_n))$
is simply the uniform law on $\Phi(x_1,\ldots, x_n)$. We deduce by comparison with Lemma \ref{L1} that for an arbitrary sequence of integers  $(x_i)_{1\leq i\leq n}$ which fulfills \eqref{E1} and $\mu(x_i)>0$, 
the distribution of the terminal state in the grabbing particle system
started with $n$ particles and such that the number of arms of the $i$-th particle is $x_i$
coincides with the conditional law $\PGW_{n,k}(\cdot \mid \Phi(x_1,\ldots, x_n))$.

To conclude the proof, we just need to recall the following well-known consequence of the ballot theorem (see, for instance, Lemma 6.1 in Pitman \cite{PiSF}).
Under the conditional law $\P_k(\cdot \mid T=n)$ of the
 Galton-Watson process started with $k$ ancestors and conditioned to have a total population with size $n$, the sequence of the number of children of individuals listed uniformly at random is given by a sequence of $n$ independent random variables with law $\mu$ conditioned to add-up to $n-k$. Recall that the event that this sequence is given by $(x_1, \ldots, x_n)$ is precisely $\Phi(x_1, \ldots, x_n)$. 
 This completes the proof of the statement. \QED

 \end{section}

\begin{section}
 {Asymptotic behavior of the empirical measure}
 
Our goal here is to investigate the asymptotic behavior of the
 empirical measure of the trees resulting from the grabbing particle system started with $n$ i.i.d. particles as $n\to\infty$. More precisely, we shall work with the following setting.

 Let $\mu$ be a critical or sub-critical  probability measure on $\Z_+$, i.e. with mean
 $$m:=\sum_{\ell=1}^{\infty}\ell \mu(\ell)\leq 1\,.$$
 We shall also always assume that $\mu(0)>0$.
Let  $\xi_1, \ldots$  be a sequence of i.i.d. variables with law $\mu$, and for
 each $n\geq 2$,  set 
 $$k(n):=n-(\xi_1+\cdots + \xi_n)\,.$$  
We work conditionally
 on the event $k(n)\geq 1$ and consider  the grabbing particle system started with $n$ particles such that the initial number of inactive arms of the $i$-th particle is $\xi_i$ for every $i\in\{1,\ldots, n\}$. We write 
 ${\bf P}_n$ for the distribution of its terminal state  in $\Phi_n:=\bigcup_{1\leq k \leq n}\Phi_{n,k}$, and 
 also  denote by $\tau_1, \ldots, \tau_{k(n)}$ the sequence  of  planar rooted trees in this terminal state, ignoring labels on vertices and edges. That is to say that the
 sequence  $(\tau_1, \ldots, \tau_{k(n)})$ is the image of the terminal state
 by the canonical projection $\Phi_{n,k(n)}\to F_{n,k(n)}$.

  We also introduce the space  of finite planar rooted trees with unlabeled edges and vertices,
  $${\mathcal T}:=\bigcup_{\ell\geq 1}F_{\ell,1}\,.$$  We denote by ${\mathbb T}^{(\mu)}$ the probability measure on ${\mathcal T}$
 which is induced by the genealogical tree of a single ancestor in a Galton-Watson process with reproduction law $\mu$.
  We are now able to state our main limit theorem.

 \begin{theorem}\label{T2}
Assume that $\mu$ is a critical or sub-critical probability measure on $\Z_+$ with $\mu(0)>0$. Then the proportion of tree structures identical to some fixed  ${\bf t}\in{\mathcal T}$ in the terminal state of the grabbing particle system converges to ${\mathbb T}^{(\mu)}({\bf t})$ in probability under ${\bf P}_n=\P(\cdot\mid k(n)\geq 1)$ as $n\to\infty$. That is, equivalently,
\begin{equation}\label{E2}
\lim_{n\to\infty}{\bf E}_n\left( \left |\frac{1}{k(n)}\sum_{\ell=1}^{k(n)} {\bf 1}_{\{\tau_{\ell}={\bf t}\}}
-{\mathbb T}^{(\mu)}({\bf t})\right | ^2\right)=0
\end{equation}

 \end{theorem}
We shall reduce the proof of \eqref{E2} to establishing the following two asymptotic behaviors under  the conditional probabilities ${\bf P}_n$.

\begin{lemma} \label{L2} The number of trees in the terminal state under the law ${\bf P}_n$ converges in probability to $\infty$, that is  
$$ \lim_{n\to\infty}{\bf P}_n(k(n)\geq k)=1\quad \hbox{  for every }k\in\N\,.$$
\end{lemma}

\proof   Note that 
$$\P(k(n)\geq k) \,=\, \P(\xi_1+\cdots + \xi_n \leq n-k)$$  tends to $1$ when  $\mu$ is  subcritical (and then the statement is obvious), but may tend to $0$ when $\mu$ is critical.  We therefore focus on the critical case.
For simplicity, we assume that the support of $\mu$ is not contained into a strict subgroup of $\Z$; the modification needed to treat the opposite case is straightforward. 
Then $\xi_1+\cdots +\xi_n-n$  is an aperiodic, irreducible recurrent random walk on $\Z$, and we know from a ratio limit theorem in Spitzer \cite{Sp} on its page 49 that for every fixed $\ell\in\Z$
$$\P(\xi_1+\cdots +\xi_n=n-\ell) \,\sim\, \P(\xi_1+\cdots +\xi_n=n)\,.$$
It follows 
$$\P(\xi_1+\cdots +\xi_n=n-\ell) \,=\, o(\P(\xi_1+\cdots +\xi_n\leq n))\,;$$ 
consequently
\begin{equation}\label{E4}
\P(\xi_1+\cdots +\xi_n\leq n-\ell) \,\sim\,\P(\xi_1+\cdots +\xi_n\leq n)\,,
\end{equation}
which yields our claim. \QED

\begin{lemma} \label{L3} We have  
\begin{equation}\label{E5}
\lim_{n\to\infty}{\bf P}_n(\tau_1={\bf t}_1, \tau_2={\bf t}_2, k(n)\geq 3)
={\mathbb T}^{(\mu)}({\bf t}_1){\mathbb T}^{(\mu)}({\bf t}_2)\,.
\end{equation}
\end{lemma}

\proof We
denote by $|{\bf t}|$ the size (number of vertices) of a tree ${\bf t}\in{\mathcal T}$.
Using Theorem \ref{T1}, we know that for each fixed $n\geq 3$
\begin{eqnarray}\label{E6}
& &{\bf P}_n(\tau_1={\bf t}_1, \tau_2={\bf t}_2, k(n)\geq 3) \nonumber\\
&=& c_n
\sum_{k=3}^{n} \P(\xi_1+\cdots + \xi_n=n-k)\PGW_{n,k}(\tau_1={\bf t}_1, \tau_2={\bf t}_2)\,,
\end{eqnarray}
with
$$c_n=1/\P(\xi_1+\cdots +\xi_n\leq n-1)\,.$$

Then recall the following celebrated identity of Dwass \cite{Dwass}:  
the first passage time process 
$$T_k=\inf\{n\geq 0: \xi_1+\cdots+\xi_n=n-k\}$$ has i.i.d. increments, which have the same distribution as the sequence
of the sizes of trees in a Galton-Watson forest with reproduction law $\mu$.
In particular the step $T_1$ is
distributed as $|\tau|$ under ${\mathbb T}^{(\mu)}$.
Further the formula of Kemperman (see, for instance,
 (6.3) in \cite{PiSF} on its page 122) states: 
 $$\P(T_k=n)=\frac{k}{n}\P(\xi_1+\cdots+\xi_n=n-k)\,,\qquad\hbox{ for every } 1\leq k < n\,.$$

By definition of the conditional law, we have thus for $k\geq 3$ and $t:=
|{\bf t}_1|+|{\bf t}_2|$
\begin{eqnarray*}
& &\PGW_{n,k}(\tau_1={\bf t}_1, \tau_2={\bf t}_2)\\
&=&{\mathbb T}^{(\mu)}({\bf t}_1)
{\mathbb T}^{(\mu)}({\bf t}_2)\frac{\P(T_{k-2}=n-t)}{\P(T_k=n)}\\
&=&{\mathbb T}^{(\mu)}({\bf t}_1)
{\mathbb T}^{(\mu)}({\bf t}_2) \frac{n(k-2)}{k(n-t)}
\frac{\P(\xi_1+\cdots + \xi_{n-t}=n-k+2-t)}{\P(\xi_1+\cdots + \xi_n=n-k)}\,,
\end{eqnarray*}
so combining with \eqref{E6}
we arrive at
\begin{eqnarray}\label{E7}
& &{\bf P}_n(\tau_1={\bf t}_1, \tau_2={\bf t}_2, k(n)\geq 3) \nonumber\\
&=& c_n {\mathbb T}^{(\mu)}({\bf t}_1)
{\mathbb T}^{(\mu)}({\bf t}_2) \frac{n}{n-t}
\sum_{k=3}^{n} \frac{k-2}{k}
\P(\xi_1+\cdots + \xi_{n-t}=n-k+2-t)\,.
\end{eqnarray}
In the subcritical case, we have $c_n\to 1$ and also, 
by the weak law of large numbers, that  the series in \eqref{E7} converges to $1$ as $n\to\infty$. This proves \eqref{E5} in that case. 

We henceforth focus on the critical case. In this direction,
from \eqref{E4} and the identity
$$\P(\xi_1+\cdots +\xi_{n+1}\leq n+1)\,=\, \sum_{k=0}^{\infty}\mu(k)\P(\xi_1+\cdots +\xi_n\leq n+1-k)$$
we get by dominated convergence that 
\begin{equation}\label{E8}
\P(\xi_1+\cdots +\xi_{n+1}\leq n+1) \,\sim\, \P(\xi_1+\cdots +\xi_n\leq n)\,.
\end{equation}
On the one hand, combining \eqref{E4}, \eqref{E8} and the obvious upper-bound
$$\sum_{k=3}^{n} \frac{k-2}{k}
\P(\xi_1+\cdots + \xi_{n-t}=n-k+2-t)\leq \P(\xi_1+\cdots + \xi_{n-t}\leq n-1-t)\,,$$
we deduce from \eqref{E7}  that
$$\limsup_{n\to\infty}{\bf P}_n(\tau_1={\bf t}_1, \tau_2={\bf t}_2, k(n)\geq 3) 
\leq {\mathbb T}^{(\mu)}({\bf t}_1)
{\mathbb T}^{(\mu)}({\bf t}_2)\,.$$
On the other hand, pick $0<\eta<1$ and let $k_{\eta}$ be an integer such that
$(k-2)/k\geq 1-\eta$ whenever $k\geq k_{\eta}$.
Combining \eqref{E4}, \eqref{E8} and the lower-bound
 $$\sum_{k=3}^{n} \frac{k-2}{k}
\P(\xi_1+\cdots + \xi_{n-t}=n-k+2-t)\geq (1-\eta)\P(\xi_1+\cdots + \xi_{n-t}\leq n-k_{\eta}+2-t)\,,$$
we deduce from \eqref{E7}  that
$$\liminf_{n\to\infty}{\bf P}_n(\tau_1={\bf t}_1, \tau_2={\bf t}_2, k(n)\geq 3) 
\geq (1-\eta){\mathbb T}^{(\mu)}({\bf t}_1)
{\mathbb T}^{(\mu)}({\bf t}_2)\,.$$
Since $\eta$ can be chosen arbitrarily small, 
this shows \eqref{E5}.
\QED

We have now all the ingredients to establish Theorem \ref{T2} using a standard argument of propagation of chaos (see e.g. Sznitman \cite{Sz}).

\proof 
 The obvious inequality
 $$\sum_{{\bf t}_1, {\bf t}_2\in{\mathcal T}}{\bf P}_n(\tau_1={\bf t}_1, \tau_2={\bf t}_2, k(n)\geq 3)={\bf P}_n(k(n)\geq 3)\leq 1 = \sum_{{\bf t}_1, {\bf t}_2\in{\mathcal T}}
 {\mathbb T}^{(\mu)}({\bf t}_1){\mathbb T}^{(\mu)}({\bf t}_2)$$
combined with Scheff\'e's lemma show that the convergence \eqref{E5}
also holds in the space $\ell^1({\mathcal T}\times {\mathcal T})$ of summable families indexed by the countable set ${\mathcal T}\times {\mathcal T}$.
It follows that 
$$\lim_{n\to\infty}{\bf P}_n(\tau_1={\bf t})
={\mathbb T}^{(\mu)}({\bf t})\,,$$
and then by an argument based on exchangeability that
$$\lim_{{n\to\infty}}{\bf E}_n\left(\frac{1}{k(n)}\sum_{\ell=1}^{k(n)} {\bf 1}_{\{\tau_{\ell}={\bf t}\}}
\right)= {\mathbb T}^{(\mu)}({\bf t})\,.$$
Using again \eqref{E5}, exchangeability and Lemma \ref{L2}, we get the asymptotic of the second moment
$$\lim_{{n\to\infty}}{\bf E}_n\left(\left(\frac{1}{k(n)}\sum_{\ell=1}^{k(n)} {\bf 1}_{\{\tau_{\ell}={\bf t}\}}\right)^2
\right)= {\mathbb T}^{(\mu)}({\bf t})^2\,,$$
which establishes \eqref{E2}. \QED

Let us briefly discuss the super-critical case, i.e. when $m>1$ (recall that we always assume that $\mu(0)>0$).
Then  $\P(k(n)\geq 1)$  tends to $0$ as $n\to\infty$, so the conditioning on this event becomes singular. More precisely, it is readily  seen  that the distribution of $k(n)$ under the conditional law ${\bf P}_n=\P(\cdot \mid k(n)\geq 1)$ converges to some geometric law. In particular Lemma \ref{L2} fails and the situation differs crucially from that describes in Theorem \ref{T2}. Nonetheless one can get an analogue of Theorem \ref{T2} in the supercritical case, provided that we condition on
$k(n)\geq cn$ for some $0<c<1$. Then standard arguments of large deviations can be 
incorporated to the  proof of Theorem \ref{T2} and yield a similar limit theorem in which $\T^{(\mu)}({\bf t})$ has to be replaced by  $\T^{(\tilde \mu)}({\bf t})$, where $\tilde \mu$
is the unique law in the exponential family of $\mu$ with mean $1-c$. Details are left to the interested reader(s).

We now conclude this work by comparing the asymptotic behaviors for the present grabbing particle system and for the model of random configuration (see \cite{BS, MR1}), starting in both cases from a large number of particles such that the sequence of the numbers of arms attached to those particles is  i.i.d. with a fixed law $\mu$.
Recall that in the latter, edges are formed by pairing arms uniformly at random, and that Molloy and Reed \cite{MR1} have shown (under slightly different assumptions) that  with high probability, the random configuration model does not possess a giant connected component if and only if $\sum_{\ell=0}^{\infty}\ell(\ell-2)\mu(\ell)\leq 0$. Further, in that case, if we select an arm uniformly at random and independently of the pairing process, then the structure of the cluster containing that arm can be described in terms of two independent Galton-Watson trees with reproduction law
$\nu(k)=(k+1)\mu(k+1)/m$; see Theorem 1 in \cite{BS} for details. So, roughly speaking,  both models produce random graphs with the same degree sequences, and for both models the structure of a typical cluster can be described in terms of Galton-Watson trees. We stress that the reproduction laws are distinct (except in the Poisson case), which may seem paradoxical as degree sequences are closely related to the latter. However one should keep in mind that selecting a cluster by
picking an arm uniformly at random differs from selecting a cluster uniformly at random,
and more precisely  induces a bias by the number of arms in the random cluster.
This explains the relation between the two reproduction laws.

 \end{section}

\vskip1cm \noindent {\bf
Acknowledgments :} We would like to thank two anonymous referees for their insightful comments on the first draft of this work. In particular, clever arguments due to a referee have permitted to remove an irrelevant condition involving stable domains of attraction in our original version of Theorem 2. 

J. B. thanks the programs
ARCUS and France-Br\'esil for supporting his visit to
IMPA at Rio de Janeiro during which the present joint work has been
undertaken.  M.E.V. acknowledges the support of CNPq, Brasil
Grant number 302796/2002-9.

\end{document}